\documentclass{article}
\usepackage{graphicx} % Required for inserting images
\usepackage{amsmath}
\usepackage{amssymb}
\usepackage{amsthm}

\usepackage{url}

\usepackage{dsfont}
\usepackage{amsfonts}
\usepackage{xspace}
\usepackage{enumitem} % Provides control over list styles

\newcommand{\zero}{\mathbf{0}}
\newcommand{\R}{\mathbb{R}}
\newcommand{\co}{\mathrm{co}}

\newcommand{\super}{\overline{\partial}}
\newcommand{\sub}{\underline{\partial}}
\newcommand{\N}{\mathbb{N}} %Neutral knots

\newcommand{\PP}{\mathbb{P}}

\newcommand{\cR}{\mathcal{R}}
\newcommand{\cQ}{\mathcal{Q}}

\newtheorem{theorem}{Theorem}[section]

\newtheorem{corollary}{Corollary}[section]

\newtheorem{remark}{Remark}[section]
\newtheorem{definition}{Definition}[section]

\title{Nonsmooth Optimisation and neural networks}
\author{Vinesha Peiris\footnote{{Centre for Optimisation and Decision Science}, {Curtin University}, {Kent Street}, {Bently}, {6102}, {Western Australia}, {Australia}, Email: m.peiris@curtin.com},  Nadezda Sukhorukova\footnote{Corresponding author, ORCID iD:  0000-0002-4078-2014, Department of Mathematics, Swinburne University of Technology, {John Street}, {Hawthorn}, {3122}, {Victoria}, {Australia}, Email: nsukhorukova@swin.edu.au}
}
\date{}

\begin{document}

\maketitle

\begin{abstract}
In this paper, we study neural networks from the point of view of nonsmooth optimisation, namely, quasidifferential calculus. We restrict ourselves to the case of uniform approximation by a neural network without hidden layers, the activation functions are restricted to continuous strictly increasing functions. We develop an algorithm for computing the approximation with one hidden layer through  a step-by-step procedure. The nonsmooth analysis techniques demonstrated their efficiency. In particular, they partially explain why the developed step-by-step procedure may run without any objective function improvement after just one step of the procedure.      
    
\end{abstract}


{\bf{MSC Classification 41A50, 26A27, 65D10, 65D12, 65D40}}

\section{Introduction}
Deep learning is a branch of machine learning which evolve rapidly within the area of artificial intelligence. Recently, the study of deep learning has been an interesting area for many researchers due to its excellent performance in a variety of practical applications~\cite{Goodfellow2016,lecun2015deep, min2017deep,raissi2019physics}, just to name a few. Deep learning focuses on using neural networks to model and solve complex problems. In general, a neural network consists of interconnected layers of nodes (input, hidden, output) that process and transform data. These networks are designed to learn patterns and representations directly from data through a process known as training, where the network adjusts its internal parameters (weights and biases) to minimise errors. These networks can also be interpreted as solving approximation problems, with their internal parameters effectively managed through optimisation techniques~\cite{Goodfellow2016,Hornik1991,SunOptDeepLearning}.

The least squares loss function (error function), widely used in neural networks, serves as a standard method to evaluate and optimise the performance of the network on a training dataset. Its popularity comes from the fact that it involves minimising a smooth quadratic objective function, making it well suited for basic optimisation techniques such as gradient descent. In some specific cases, the loss function can be adapted to a uniform (Chebyshev) approximation based model, which minimises the maximum error and outperforms least squares--based models when the training dataset is accurate, but limited in size~\cite{ peiris2022rational,PeirisRoshSuk2024}. 

In this research, we use a simple model of an artificial neural network that has no hidden layers and only one output node. In the case of deep neural networks, his model corresponds to the optimisation between two consecutive  layers.  We explore the use of two types of activation functions for the output layer: Leaky ReLU function and strictly increasing smooth functions. We utilise the uniform approximation based loss function which leads to a nonsmooth and nonconvex optimisation problem and therefore, nonsmooth optimisation techniques are needed to address them. One of the most comprehensive textbooks on nonsmooth optimisation is~\cite{bagirov2014introduction}. 

The quasidifferential is an appropriate tool for the construction of the solution of nonsmooth and nonconvex optimisation problems. It is a possible extension of the notion of gradient to the case of nonsmooth functions. Since its introduction by V.F. Demyanov, A. M. Rubinov and L. N. Polyakova~\cite{dempol1980minimization,  dempolalex1986nonsmoothness,demrub1980quasidifferentiable,demyanov1995introduction}, it has been an interesting research direction for many mathematicians. By using the well-developed calculus rules for quasidifferential functions~\cite{DemyanovDixon,Demyanov2000}, we give a characterisation of necessary optimality condition for our optimisation problems which proves that every local minimum is an inf-stationary point. Moreover, we propose applying the bisection algorithm from~\cite{PeirisRoshSuk2024} step-by-step to solve the corresponding optimisation problem.

This paper is organised as follows. The motivation for the current study is discussed in Section~\ref{sec:motivation}. Section~\ref{sec:quasidifferentials} explains the notion of quasidifferentials and related calculus rules. In Section~\ref{sec:math_formulation}, we introduce the mathematical formulation of the networks that we use. In Section~\ref{subsec:step-by-step}, we discuss the step-by-step process to solve the optimisation problems and we provide the necessary condition for optimality. Section~\ref{sec:experiments} provides numerical experiments and finally, Section~\ref{sec:conclusion} explains conclusions and future research directions.

\section{Motivation} \label{sec:motivation}

%\subsection{Neural network and its mathematical background}

The theory of  Artificial Neural Networks (ANNs) is based on solid mathematical background.    The origins come back to the work of A.~Kolmogorov and his student V.~Arnold. Their   celebrated Kolmogorov-Arnold representation Theorem~\cite{Arnold57,Kol57} is an attempt to solve the 13th problem of Hilbert. Then, further development of this result lead to a convenient and efficient adaptation to function and data approximation, established in~\cite{Cybenko,Hornik1991,LeshnoPinkus1993,pinkus_1999}. An excellent overview of optimisation techniques for ANNs can be found in~\cite{SunOptDeepLearning}. 

Essentially, ANNs construct approximation in the form of a sum
%\begin{equation}\label{eq:sum}
$$
\sum_{i=1}^{n}a_i\sigma((b_1,b_2,\dots,b_m){\bf T}+b_0),   
%\end{equation}
$$
where ${\bf T\in \R^m}$ and $$a_1,\dots,a_n,b_0,\dots, b_m$$
are the parameters of approximation and these parameters are the decision variable of the optimisation problems. The goal is to minimise the loss functions, which are the approximation error (least squares, maximal absolute deviation, sum of absolute deviations, etc.). 

Most ANN packages rely on least squares, since the corresponding optimisation problems are smooth and can be locally optimised. In our research, we optimise the maximal deviation (uniform norm) and the corresponding optimisation problems are nonsmooth and nonconvex. 

ANN are very popular techniques in modern numerical applications, in particular, for construction accurate approximation~\cite{approxdiffeq,variational, approxKorobov,PeirisRoshSuk2024} just to name a few. Still there are a number of open problems and avenues for improvements~\cite{limitation}.

Due to the nonsmoothness of the problem, one needs nonsmooth optimisation techniques to tackle them. One possible way is to do it using quasidifferential calculus~\cite{demyanov1995introduction}. This is just one of the possible extension of the notion of gradient to the case of nonsmooth functions. Quasidifferential calculus has been applied to a number of difficult nonsmooth optimisation problems, including best  free knot polynomial spline  approximation~\cite{SukUgonTrans2017}.  

%\section{Nonsmooth optimisaton}
\section{Quasidifferentiable calculus} \label{sec:quasidifferentials}

Nonsmooth optimisation is a branch of optimisation where one deals with nonsmooth functions. The history of nonsmooth optimisation can be traced back to P.~Chebyshev and his approximation problems~\cite{Cheb}. A very detailed review of the development on this field can be found in~\cite{Gaudioso2022}. In this study, we are using the notion of quasidifferentiability, one of the  possible generalisation of smoothness, originally described in~\cite{dempol1980minimization,demrub1980quasidifferentiable} and later developed to a rich and fruitful field~\cite{DemyanovDixon,demyanov1995introduction,Demyanov2000}. This theory allows one to compute descent directions for nonsmooth functions. If no descent direction can be obtained, the point is {\it inf-stationary in the sense of Demyanov-Rubinov}. The formal definition is provided later in this section. This can be done through the so-called Quasidifferential Calculus, developed in~\cite{demyanov1995introduction}.

In this paper, we only present the essential results, directly related to our study. For more information, one should refer to \cite{demyanov1995introduction,Demyanov2000}.

          \begin{definition}[{\cite[Chapter 1, page 7]{Demyanov2000}}]
            \label{def:quasidifferential}
            A function \(\Psi\) defined on an open set \(\Omega\) is
            \emph{quasidifferentiable} at
            a point \(X \in \Omega\) if it is locally Lipschitz continuous, directionally
            differentiable at this point and there exists compact, convex sets
            \(\sub \Psi(X)\) and \(\super \Psi(X)\) such that %for any direction \(g\) the
            the derivative of \(\Psi\) at \(X\) in any direction \(g\) can be expressed as
            \[
      	  \Psi'(X,g)=\max_{\mu \in \sub{\Psi(X)}} \langle \mu,g \rangle +
              \min_{\nu \in \super{\Psi(X)}} \langle \nu,g \rangle.
            \]
            The sets \(\sub \Psi(X)\) and \(\super \Psi(X)\) are called respectively the
            \emph{sub-}
            and \emph{superdifferential} of the function \(\Psi\) at the point \(X\). The pair \([\sub
      	  \Psi(X),\super \Psi(X)]\) is called a \emph{quasidifferential} of the function \(\Psi\) at
            the point \(X\).
          \end{definition}

      A point \(X^*\) satisfying condition~\eqref{eq:quasidiffcondition} is an
      \emph{inf-stationary} point in the sense of Demyanov-Rubinov~\cite{demyanov1995introduction}:   \begin{equation}\label{eq:quasidiffcondition}
        -\overline{\partial}\Psi(X^*)\subset \underline{\partial}\Psi(X^*).
      \end{equation}
       At any local minimiser \(X^* \in \Omega\)
      of a quasidifferentiable function \(\Psi\) the condition~(\ref{eq:quasidiffcondition}) is satisfied (necessary condition for local minimality).
From the definition, one can see that  there may be several ways to choose quasidifferentials. For example, in the case of smooth functions, the subdifferential can be chosen as the gradient while the superdifferential is zero vector of the corresponding dimension. Another possibility is to choose the gradient as the supperdifferential and the zero vector of the corresponding dimension as the subdifferential.

For this study, we will need the expression for calculating the quesiifferential of the sum and maximum of quasidifferential function and also the quasidifferential of the product of a 
 quasidifferential function and a real number. For more details, one can refer to~\cite{demyanov1995introduction}.
\begin{itemize}
    \item{\bf Property 1:} Let $f_1(x)$ and $f_2(x)$ be two quasidifferential functions at $x$  whose  quasidifferentials are $[\underline{\partial}f_1(x),\overline{\partial}f_1(x)]$ and  $[\underline{\partial}f_2(x),\overline{\partial}f_2(x)]$, then the quasidifferential of $f_1(x)+f_2(x)$ is    $$[\underline{\partial}f_1+\underline{\partial} f_2(x),\overline{\partial}f_1(x)+\overline{\partial}f_2(x)].$$
    \item{\bf Property 2:}   Let $f(x)$ be a quasidifferential function  at $x$, whose quasidifferential is $[\underline{\partial}f(x),\overline{\partial}f(x)]$ and a real number $\lambda$, then the quasidifferential of the function $\lambda f$ is
        \begin{itemize}[label=\textbullet]
            \item $[\lambda\underline{\partial}f(x),\lambda\overline{\partial}f(x)]$ if $\lambda\geq 0$;
             \item $[\lambda\overline{\partial}f(x),\lambda\underline{\partial}f(x)]$ if $\lambda< 0$.
        \end{itemize}
    \item{\bf Property 3:} Let $f_1,\dots,f_m$ defined on an open set~$X$ are quasidifferential in $x\in X$ and 
    $$\phi=\max_{i=1,\dots, m}f_i(x)$$
    then $\phi$ is also quasidifferential at $x$ and 
    $$\underline{\partial}\phi(x)=
 \co\cup_{k\in R(x)}\left (\underline{\partial}f_k(x)-\sum_{i\in R(x),~i\neq k}\overline{\partial}f_i(x)\right ),$$
 $$\overline{\partial}\phi(x)=
\sum_{k\in R(x)}\overline{\partial}f_k(x),$$
    where $R(x)=\{i\in I|f_i(x)=\phi(x)\}$ and $I=\{1,\dots,m\}$. Essentially, the set of indices~$R(x)$ corresponds to the indices, where the maximum is reached. 
    \item{\bf Property 4:} This property originally appeared in~\cite{minimax}, p.~243 (in Russian). It can be also found in~\cite{demyanov1995introduction}.  Let $X$ be an open set in~$\R^n$ and $Y\subset \R^m$ is a compact set. Function $f(x,y)$ is defined on $X\times Y$ and continuous on this set together with the partial derivative $f'_x(x,y)$. Let $$\phi(x)=\max_{y\in Y}f(x,y),$$ then a quasidifferential can be constructed as 
    $$\underline{\partial}\phi(x)=\co\{v\in\R^n| v=f'_x(x,y),~y\in R(x)\},$$
    where $R(x)=\{y\in Y|f(x,y)=\phi(x)\}$
    and $\overline{\partial}\phi(x)$ is a zero vector of the corresponding dimension.
\end{itemize}

\section{Mathematical formulation} \label{sec:math_formulation}
%\subsection{Mathematical formulation}
This problem can be formulated as a continuous problem (function approximation) or discrete probelm (data approximation). The latter is more common in practical problems, since the exact function expression is unknown. This is especially clear in classification problems, where each entry point belongs to a certain class (label). 

During the training phase, one constructs a network, whose weight are optimised to match the data and the corresponding label (class). Then, during the test phase, one uses the constructed network (training phase) to predict the class for the data points that were not used in the training phase.  Essentially, in a dataset, every entry is a vector (features) and the function value is the class.

The discrete version of the problem is as follows:
% \begin{equation}\label{eq:free_disc}
%  {\text{minimise}}  \max_{\substack{{\bf T}_j\\  j=1,\dots,N}}  \left| \sum_{i=1}^{n}a_i\sigma({\bf w}^i{\bf T}_j+w_{0}^i)- f({\bf T}_j)\right|, \end{equation}
%  \begin{equation}\label{eq:free_disc_c0nstraints}
%   {\text{subject to}}~ X\in \R^{n+nd+n},
% \end{equation}
% where $N$ is the number of discretisation points (number of dataset entries), ${\bf T}_j\in Q$, $j=1,\dots,N$ are $d$-dimensional vectors (dataset points), $Q$ is a $d$-dimensional hypercube,  $n$ is the maximal number of affine pieces, weights~$w_k^i$, $i=0\dots,n$, $k=1,\dots,d$, vectors~${\bf w}^i=(w_1^i,\dots,w_d^i)^T$ and coefficients~$a_i$, $i=1,\dots,n$ are the decision variables. For convenience, all the decision variables are placed in a single vector 
% $X=(a_1,\dots,a_n, w_{0}^1,\dots,w_{d}^{n})$. $\sigma(x)$ is the activation function and $f$ is the original function (function to approximate). 
\begin{equation}\label{eq:free_disc}
 {\text{minimise}}  \max_{\substack{{\bf T}_j\\  j=1,\dots,N}}  \left| \sum_{i=1}^{n}a_i\sigma({\bf w}^i{\bf T}_j+w_{0}^i)- f({\bf T}_j)\right|, \end{equation}
 \begin{equation}\label{eq:free_disc_c0nstraints}
  {\text{subject to}}~ X\in \R^{n+nd+n},
\end{equation}
where 
\begin{itemize}
    \item $N$ is the number of discretisation points (number of dataset entries),
    \item ${\bf T}_j\in Q$, $j=1,\dots,N$ are $d$-dimensional vectors (dataset points) where $Q$ is a $d$-dimensional hypercube,
    \item $n$ is the maximal number of affine pieces,
    \item vectors~${\bf w}^i=(w_1^i,\dots,w_d^i)^T$ are the weights where $w_k^i$, $i=0,\dots,n$, $k=1,\dots,d$, and coefficients~$a_i$, $i=1,\dots,n$ are the decision variables. For convenience, all the decision variables are placed in a single vector, 
    $$X=(a_1, \dots, a_n, w_{1}^1, \dots, w_{d}^{n}, w_{0}^1, \dots, w_{0}^n),$$
    \item $\sigma(x)$ is the activation function and, 
    \item $f$ is the original function (function to approximate).
\end{itemize}

In this paper, we use two types of activation function: 
\begin{itemize}
\item
Leaky ReLu $\sigma(x)=\max\{x,\alpha x\}$, where $\alpha$ is typically between 0.01 and 0.03.
\item Strictly increasing smooth function, for example, a sigmoidal function.
\end{itemize}

Optimisation problem~(\ref{eq:free_disc})-(\ref{eq:free_disc_c0nstraints}) is nonsmooth and non-convex and therefore it is a challenging problem. In the case when there is only one affine piece ($n=1$) and the activation function is strictly monotone (increasing functions are more common in practical problems) the problem is quasiconvex and can be efficiently solved using the so called bisection method~\cite{PeirisRoshSuk2024}. {\bf Moreover, this method leads to a global minimum}.

In this paper, we will continue the study of simple neural networks without any hidden layers. The input (data point) goes directly to the output. In particular, we will find the affine pieces step-by-step. 

\section{Step-by-step algorithm} \label{subsec:step-by-step}

In this paper we propose a simple procedure of applying the bisection algorithm  from~\cite{PeirisRoshSuk2024} step-by-step: first we approximate the original function and then the deviation function. One possible extension is to determine the affine coefficients, fix them and refine the coefficients~$a_i$ rather than fixing them at~$1$.   
\begin{enumerate}
\item Solve (\ref{eq:free_disc})-(\ref{eq:free_disc_c0nstraints}) for $n=1$ (quasiconvex problem, see~\cite{PeirisRoshSuk2024} for details). 
\item Subtract the obtained approximation from the original function and this is our new function for approximation. 
\end{enumerate}

In the next section, we study the theoretical properties of inf-stationary points of~(\ref{eq:free_disc}). This study is essential for understanding the properties of the proposed step-by-step algorithm.

\subsection{Inf-stationary points and their properties}

In this section, we refer the the following optimisation problem, which is a special case of ~(\ref{eq:free_disc})-(\ref{eq:free_disc_c0nstraints}) when $n=1$:

\begin{equation}\label{eq:free_disc_simple}
 {\text{minimise}}  \max_{\substack{{\bf T}_j\\  j=1,\dots,N}}  \left| a\sigma({\bf w}{\bf T}_j+w_{0})- f({\bf T}_j)\right|, \end{equation}
 \begin{equation}\label{eq:free_disc_c0nstraints_simple}
  {\text{subject to}}~ X\in \R^{d+2},
\end{equation}

Note that function
$$\max_{\substack{{\bf T}_j\\  j=1,\dots,N}}  \left| a\sigma({\bf w}{\bf T}_j+w_{0})- f({\bf T}_j)\right|$$
can be rewritten as follows:
\begin{equation}\label{eq:free_disc_simple_difference}
\max_{\substack{{\bf T}_j\\  j=1,\dots,N}}  \left\{ a\sigma({\bf w}{\bf T}_j+w_{0})- f({\bf T}_j),f({\bf T}_j)-a\sigma({\bf w}{\bf T}_j+w_{0})\right\}.
\end{equation}

In the rest of the section, give a characterisation of inf-stationary points of function from~(\ref{eq:free_disc_simple}), where $a=1$. 

Note that {\bf every local minimum is an inf-stationary point} (necessary condition)~\cite{demyanov1995introduction}.

The case when $a$ is part of decision variables is out of scope of this paper, but will be considered in our future research.

Let ${\bf w}=(w_1,\dots,w_d)$ and 
$${\bf W}=\left[ {\begin{array}{l}
    w_0  \\
    {\bf w}  \\
  \end{array} } \right].$$

\begin{definition}
Function $\delta({\bf W}, {\bf T}_j)=\sigma({\bf w}{\bf T}_j+w_{0})- f({\bf T})$ is called the {\it deviation function}, while  $\Delta({\bf W})=\max_{{\bf T}_j\in D} |\sigma({\bf w}{\bf T}_j+w_{0})- f({\bf T})|$.
\end{definition}
\subsection{Smooth activation function}

Due to Property~4, inf-stationary points of~(\ref{eq:free_disc_simple_difference}) are as follows:
\begin{equation}
    \zero_{d+1}\in \co\{v\in\R^n| v=\delta'_{\bf W}({\bf W}, {\bf T}),~{\bf T}\in R({\bf W})\},
    \end{equation}
where
%$${\bf W}=\left[ {\begin{array}{l}
%    w_0  \\
%    {\bf w}  \\
%  \end{array} } \right]$$
$\cR({\bf W})=\{{\bf T} \in D | \max_{T\in D}|\delta({\bf W}, {\bf T})|=\Delta\}$, here $D$ is the collection of discretisation points $T_i$, $i=1,\dots,N$ (dataset). Note that in this problem, the set~$R({\bf W})$ corresponds to the set, where the maximum of absolute deviation is reached. The deviation itself can be positive or negative and it corresponds to the sign of $\delta({\bf W}, {\bf T})$. Properties~1 and~3 are also used in computing quasidifferentials. Therefore, a necessary optimality condition (inf-stationary points) can be rewritten as follows:
\begin{equation}\label{eq:noc}
    \zero_{d+1}\in \co\{v\in\R^n| v=\delta'_{\bf w}({\bf W}, {\bf T}),~{\bf T}
    \in (\N({\bf W})\cup \PP({\bf W}))\},
    \end{equation}
    where $\N({\bf W})$ is the set of maximal absolute deviation with the negative sign ($\delta({\bf W}, {\bf T})<0$), while $\PP({\bf W})$ is the set of maximal absolute deviation with the positive sign ($\delta({\bf W}, {\bf T})\geq 0$). Let $N$ be the total number of data points in the dataset $D$, that is, $D=\{{\bf T}_j, j=1,\dots, N\}$. Then condition~(\ref{eq:noc}) is equivalent to the following:
    \begin{equation}\label{eq:noc1}
\zero_{d+1}\in \co\{A,B\},
\end{equation}
where
$$A=\{v\in\R^n| v={\partial \sigma\over \partial({\bf w}{\bf T}_j+w_{0})}{\bf E}_j),~{\bf T_j}\in \PP({\bf W})\},$$
$$B=\{v\in\R^n| v=-{\partial \sigma\over \partial({\bf w}{\bf T}_j+w_{0})}{\bf E}_j),~{\bf T_j}\in \N({\bf W})\},$$
\begin{equation}\label{eq:Ej}
{\bf E}_j=\left[ {\begin{array}{l}
    1  \\
    {\bf T}_j  \\
  \end{array} } \right].
\end{equation}
    
%\zero_{d+1}\in \co\{v\in\R^n| v={\partial \sigma\over \partial({\bf w}{\bf T}_j+w_{0})}{\bf E}_j),~{\bf T_j}\in \PP({\bf T})
%    \in (\N({\bf T})\cup \PP({\bf T}))\},
    
 %   \zero_{d+1}\in \co\{(v\in\R^n| v=\delta'_{\bf w}({\bf w}, {\bf T}),~{\bf T}
  %  \in (\N({\bf T}) \cup(v\in\R^n| v=\delta'_{\bf w}({\bf w}, {\bf T}),~{\bf T}
   % \in  \PP({\bf T}))
   % \in \PP(x))\}.
    %\end{equation*}
    
    Let~$N_1$ be the total number of maximal deviation points for a given set of weights. Assume that the maximal deviation points are re-numerated in such a way that ${\bf T}_j$, $j=1,\dots,N_2$ correspond to positive deviation, while ${\bf T}_j$, $j=N_2+1,\dots,N_1$ correspond to negative deviation. Then condition~(\ref{eq:noc1}) is equivalent to the existence of non-negative coefficients $\alpha_i\geq 0$, $i=1,\dots,N_1$, such that at least one of them is strictly positive and 
    \begin{equation}\label{eq:noc_original}    \zero_{d+1}=\sum_{i=1}^{N_2}{\partial \sigma\over \partial({\bf w}{\bf T}_j+w_{0})}{\bf E}_i
    -\sum_{i=N_2+1}^{N_1}{\partial \sigma\over \partial({\bf w}{\bf T}_j+w_{0})}{\bf E}_i.
    \end{equation}

In the case of {\bf smooth activation function}, regardless of their monotonicity, inf-stationary points in the sense of Demyanov-Rubinov are described in Theorem~\ref{thm:smooth}.

\begin{theorem}\label{thm:smooth}
Assume that the activation function is smooth. 
Point ${\bf W}$ is inf-stationary in the sense of Demyanov-Rubinov if and only if the sets 
$$\co\left\{{\partial \sigma\over \partial({\bf w}{\bf T}_j+w_{0})}{\bf E}_i,~i=1,\dots,N_2\right\}$$
and
$$\co\left\{{\partial \sigma\over \partial({\bf w}{\bf T}_j+w_{0})}{\bf E}_i,~i=N_2+1,\dots,N_1\right\}$$
intersect.
\end{theorem}
{\bf Proof:} 
Consider the first coordinate of ${\bf E}_j$, $j=1,\dots,N_1$, it is clear that
$$\sum_{i=1}^{N_2}\alpha_i=\sum_{j=N_2+1}^{N_1}\alpha_j=K>0.$$

Therefore,
$$\sum_{i=1}^{N_2}{\alpha_i\over K}=\sum_{j=N_2+1}^{N_1}{\alpha_j\over K}=1.$$
Hence, there exists point 
$${\bf v}=\sum_{i=1}^{N_2}{\alpha_i\over K}{\partial \sigma\over \partial({\bf w}{\bf T}_j+w_{0})}{\bf E}_i=\sum_{j=N_2+1}^{N_1}{\alpha_j\over K}{\partial \sigma\over \partial({\bf w}{\bf T}_j+w_{0})}{\bf E}_j.$$

Therefore, sets $$\co\left\{{\partial \sigma\over \partial({\bf w}{\bf T}_j+w_{0})}{\bf E}_i,~i=1,\dots,N_2\right\}$$
and
$$\co\left\{{\partial \sigma\over \partial({\bf w}{\bf T}_j+w_{0})}{\bf E}_i,~i=N_2+1,\dots,N_1\right\}$$
intersect and this condition is equivalent to~(\ref{eq:noc_original}).
\qed
\begin{remark}
Theorem~\ref{thm:smooth} is an important theoretical result, which can be also seen as a multidimentional generalisation of the notion of the alternating sequence. In~\cite{Sukhorukova2018MATRIX} a similar generalisation was developed for multivariate linear approximation. In practice, however, these conditions are not very easy to verify. In particular, this due to the fact that it is hard to develop an efficient procedure for finding  all the maximal deviation points, especially in multivariate settings. 
\end{remark}
\subsection{Leaky ReLU}
Note, first of all, that Leaky ReLU function is positively homogeneous. Therefore, if $a$ is non-negative, this decision variable is redundant, since it can be assigned to~1 and the weights $w_0,\dots,w_d$ will be readjusted. Essentially, one needs to consider only two possibilities: $a$ is non-negative and $a$ is negative. Similar to the case of smooth activation functions, we assume that $a=1$, but leave other options for our further research directions.

 In the case of Leaky ReLU as the activation function, the optimisation problem is as follows.
\begin{equation}\label{eq:free_Leaky}
\min_{\bf W}\phi({\bf W})=\max_{\substack{{\bf T}_j\\  j=1,\dots,N}} 
\left\{ 
H_j({\bf W}),-H_j({\bf W})
 \right\},
\end{equation}
where
$$
H_j({\bf W})=\max\{
\alpha({\bf w}{\bf T}_j+w_{0}),{\bf w}{\bf T}_j+w_{0}\}- f({\bf T}_j).
$$
It is clear that functions $H_j({\bf W})$ is non-smooth and therefore Theorem~\ref{thm:smooth} is not applicable. Hence, we suggest another approach, which is based on the direct calculation of the quasidifferential calculus developed in~\cite{demyanov1995introduction}. In this section we use Properties~1-3.

Let $\Delta({\bf W})$ be the maximal absolute deviation for the set of weights~${\bf W}$:  
$$\Delta({\bf W})=\max_{{\bf T}_j,  j=1,\dots,N} 
\left\{ 
H_j({\bf W}),-H_j({\bf W})
 \right\}.$$
Consider only the indices where the maximum is reached: 
\begin{equation}
    \cR({\bf W})=\{j\in 1,\dots,N|H_j({\bf W})=\Delta({\bf W})\}
\end{equation}
and 
\begin{equation}
    \cQ({\bf W})=\{j\in 1,\dots,N|-H_j({\bf W})=-\Delta({\bf W})\}.
\end{equation}
The quasidifferential can be constructed as follows (using Properties~1-3):
\begin{equation}\label{eq:superdiff}
    \overline{\partial}\phi({\bf W})=
    \sum_{j\in \cQ({\bf W})}\overline{\partial}(-H_j({\bf W}))=-\sum_{j\in \cQ({\bf W})}\co\{\alpha {\bf E}_j,{\bf E}_j\};
\end{equation}
\begin{equation}\label{q:subdiff}
    \underline{\partial}\phi({\bf W})=
    \co\bigcup_{i\in\cR({\bf W})}\left(
    \co\{\alpha {\bf E}_i,{\bf E}_i\}-\left(\zero+\sum_{j\in\cQ({\bf W}),~j\neq i}\co\{\alpha{\bf E}_j,{\bf E}_j\}\right)
    \right),
\end{equation}
where ${\bf E}_i$, $i\in \cR({\bf W})\cup\cQ({\bf W})$ is the same as in~(\ref{eq:Ej}). Since the sets $\cR({\bf W})$ and $\cQ({\bf W})$ are disjoint, the summation indices of the  subdifferential can be simplified:
\begin{equation}\label{q:subdiffsimple}
   \underline{\partial}\phi({\bf W})=
    \co\bigcup_{i\in\cR({\bf W})}\left(
    \co\{\alpha {\bf E}_i,{\bf E}_i\}-\left(\zero+\sum_{j\in\cQ({\bf W})}\co\{\alpha{\bf E}_j,{\bf E}_j\}\right)
    \right).
\end{equation}
 From the explicit calculation of the quasidifferential, one can see that the location of the maximal deviation points determine whether the chosen weights correspond to an inf-stationary point or not. 

 \begin{theorem}\label{thm:leakyReLU}
     Let $f_1({\bf T})$ and $f_2({\bf T})$ are two distinct functions and the corresponding weight vectors are ${\bf W}_1$ and ${\bf W}_2$. If the sets of maximal deviation points for these two approximations coincide (the sign of deviation is relevant) then ${\bf W}_1$ is a stationary point if and only if ${\bf W}_2$ is a stationary point. 
 \end{theorem}
 {\bf Proof:} 
 Assume that ${\bf W}_1$ is an inf-stationary point in the sense of Demyanov-Rubinov, then the inclusion~(\ref{eq:quasidiffcondition}) is satisfied. Since the maximal deviation points are the same for both approximations, the corresponding vectors ${\bf E}$, coincide and therefore the inclusion~(\ref{eq:quasidiffcondition}) is also satisfied for ${\bf W}_2$. The converse can be proven in the same manner.
 \qed
 
 \begin{corollary}\label{cor:stepbystep}
  In the case  of Leaky ReLU activation function, step-by-step bisection terminates at an inf-stationary point at each step of the procedure.  
 \end{corollary}

 Corollary~\ref{cor:stepbystep} is very important for the development of any step-by step procedure, which may appear as a reasonable and straightforward approach. Namely, since each step of this procedure terminates at an inf-stationary point of the objective function, any additional step of the procedure may be irrelevant, since it may move from one solution to another without the improvement of the objective function. The objective function is quasiconvex and therefore we can not guarantee that any in-stationary point is a global minimiser, but the bisection procedure always terminates at a point, where the values of the objective function is within $\varepsilon$  of the global minimiser.

\section{Numerical experiments} \label{sec:experiments}

In this section we present the results of two numerical experiments with the proposed step-by-step procedure. In both cases, we work on optimisation problem~(\ref{eq:free_disc})-(\ref{eq:free_disc_c0nstraints}), where $n=1$, ${\bf T}_j$, $j=1,\dots,N$ are the points  from~\cite{dataset}. The values of the function~$f$ is the class label~$\{0,1\}$. We use two datasets:  TwoLead-ECG test (82 features, 1139 points) and TwoLead-ECG training (82 features, 82 points). In both experiments, we use parameter~$\alpha=0.01$ for Leaky ReLU and the bisection tolerance~$\varepsilon=10^{-6}$.

Note that the step-by-step procedure (for $n$ steps) is not equivalent to the direct implementation of optimisation of~(\ref{eq:free_disc})-(\ref{eq:free_disc_c0nstraints}) for $n>1$. The direct optimisation applied to this problem is hard even for modern optimisation techniques and the step-by-step procedure is considered as a suitable compromise between the accuracy of approximation and the efficiency of the corresponding optimisation methods. 
\subsection{TwoLead-ECG test}
The first step of the step-by-step procedure terminates at the point with the objective function value~0.3163. The second and following application of the step-by-step procedure leads to a different point without improving the values of the objective function (within the predefined tolerance~$\varepsilon$). This result can be explained by Corollary~\ref{cor:stepbystep}, since the starting point for the second step is an inf-stationary point.

\subsection{TwoLead-ECG training}

This dataset is smaller and the first step of the bisection procedure leads to interpolation (with the tolerance $\varepsilon=10^{-6}$). The maximal absolute deviation is~$9.4414*10^{-7}$, which is less than the tolerance. 

The second step of the procedure improves the maximal absolute deviation to~$9.3942*10^{-7}$.  There is no contradiction with Corollary~\ref{cor:stepbystep}, since the starting point does not have to be inf-stationary. The application of the same procedure again leads to the maximal absolute deviation of~$1.4877*10^{-14}$, which is lower than the tolerance~$\varepsilon=10^{-6}$, used in the bisection method.

\section{Conclusions and future research directions} \label{sec:conclusion}

In this paper, we study artificial neural networks from the point of view of nonsmooth optimisation and quasidifferential calculus, developed in~\cite{demyanov1995introduction}. We develop necessary optimality conditions that are based on the notion of inf-stationary points. 

We also propose a step-by step algorithm, where at each step the bisection method is used. This method leads to a global minimum when only one affine piece is used. In order to add one more affine piece, the procedure is applied for approximating the deviation function by the next piece. In some of our numerical experiments, the algorithm found more than one piece, while is some others it could not escape from the starting point. This can be partially explained by the fact that the starting point is inf-stationary.

While the current paper is mostly theoretical, in the future we will work on the algorithm and the improvement of its efficiency. In particular, it may be possible to escape from an inf-stationary point by changing the activation function, in particular to cover the case when the coefficient $a$ is negative.

\section*{Acknowledgments}
We are grateful to the Australian Research Council for supporting this work via Discovery Project DP180100602.

Data sharing is not applicable to this article as no data sets were generated or analysed during the current study.
\appendix
\section{Appendix}
In the Appendix we present the essential information on the application of a Bisection method to Artificial neural networks used in~\cite{PeirisRoshSuk2024}. The optimisation problem is as follows
\begin{equation*}
\min_{w \in \R^{n+1}}\max_{i\in 1:N} \left|y^i - \sigma\left(\sum_{j=1}^n w_j x_j^i+w_0\right)\right|.
\end{equation*}
Here $ w = (w_0,w_1,\dots, w_n) = (w_0,\bar w)\in \R\times \R^n$ are the weights to be decided, and $Z = \{(\bar x_i, \bar y_i)\}_{i=1}^N$, with  $(\bar x^i,\bar y^i)\in \R^n\times \R$, $i\in \{1,\dots, N\}$ is the training set. 

We know  that the max of quasiconvex functions  is quasiconvex, hence we can apply the bisection method. 

The objective function is nonnegative (absolute deviation), so we can choose $l_0= 0$ as the lower bound for the optimal value. For the upper bound we can substitute any value of the parameter $w$ in the objective, for instance, $w = 0\in \R^{n+1}$, then  
\begin{equation*}
u_0 := \max_{i\in 1:N} \left|\bar y^i- \sigma(0)\right|.
\end{equation*}

%We know that the optimal value of the objective function is between $l_0$ and $u_0$. 
Let $L_1:= \frac{l_0+u_0}{2}$. On each iteration of the algorithm, we solve the feasibility problem 
\begin{equation}\label{eq:feasibility}
\max_{i\in 1:N} \left|\bar y^i - \sigma\left(\sum_{j=1}^n w_j \bar x_j^i+w_0\right)\right|\leq L_k.
\end{equation}
If the problem is feasible, we assign $l_k = l_{k-1}$, $u_{k} = L_k$. Otherwise we assign $l_k = L_{k}$, $u_k = u_{k-1}$. 

We set a threshold $\varepsilon>0$.  The stopping criterion is 
$$u_{k_0}-l_{k_0}<\varepsilon.$$ 

This procedure  terminates in a finite number of steps. The convergence is linear~\cite{SL}.

The feasibility problem \eqref{eq:feasibility} %that we are required to solve on every iteration 
can be equivalently rewritten as 
\begin{equation*}
\bar y^i - L_k \leq  \sigma\left(\sum_{j=1}^n w_j \bar x_j^i+w_0\right) \leq \bar y^i + L_k \quad \forall i \in \{1,\dots, N\}.
\end{equation*}
Since $\sigma$ is strictly increasing and therefore invertible, this can be rewritten as 
\begin{equation*}
\sigma^{-1}(\bar y^i - L_k) \leq  \sum_{j=1}^n w_j \bar x_j^i+w_0 \leq \sigma^{-1}(\bar y^i + L_k) \quad \forall i \in \{1,\dots, N\}.
\end{equation*}
%Note that if $\varphi$ is nondecreasing, we can replace the inverses with the relevant minima and maxima over the value of the multifunction inverse, as in
%\[
%\inf(\varphi^{-1}(y^i - L_k)) \leq  \sum_{j=1}^n w_j^T x_j^i+w_0 \leq \sup(\varphi^{-1}(y^i + L_k)). 
%\]
Therefore, at each step of the algorithm, one needs to solve linear programming problems and therefore it can be solved efficiently with any standard linear programming technique.

%\bibliographystyle{plain}
%\bibliography{references.bib}

\end{document}